# An integrated rolling horizon approach to increase operating theatre efficiency

Belinda Spratt* and Erhan Kozan

*Queensland University of Technology, 2 George St, Brisbane, QLD 4000, Australia*
*E-mail: b.spratt@qut.edu.au; e.kozan@qut.edu.au*

**Abstract**

Demand for healthcare is increasing rapidly. To meet demand, we must improve the efficiency of our public health services. We present a mixed integer programming (**MIP**) formulation that simultaneously tackles the integrated Master Surgical Schedule (**MSS**) and Surgical Case Assignment (**SCA**) problems. We consider volatile surgical durations and non-elective arrivals whilst applying a rolling horizon approach to adjust the schedule after cancellations, equipment failure, or new arrivals on the waiting list. A case study of an Australian public hospital with a large surgical department is the basis for the model. The formulation includes significant detail and provides practitioners with a globally implementable model. We produce good feasible solutions in short amounts of computational time with a constructive heuristic and two hyper metaheuristics. Using a rolling horizon schedule increases patient throughput and can help reduce waiting lists.

## 1. Introduction

The efficient utilisation of hospital resources is becoming increasingly important in both private and public hospitals. An increase in patient demand due to aging populations, improved screening techniques and wider access to medical care has led to long waiting times for public patients. Operations Research has been applied to a wide range of healthcare problems including planning, logistics, and practice (Rais and Viana, 2011). In this paper, we focus on the efficient scheduling of the surgical department in public hospitals. Poor planning of the surgical department has a major impact on downstream wards and can reduce the efficiency of the entire hospital.

Operations research techniques improve the scheduling of surgical procedures with respect to a wide variety of performance measures. For a detailed review on operating theatre (**OT**) planning and scheduling, including categorisation of articles by patient characteristics, performance measures, and decision delineation (among other fields) see Samudra et al. (2016) and Zhu et al. (2018). Samudra et al. (2016) place particular emphasis on the current challenges and pitfalls observed in recent literature. Zhu et al. (2018) show that most of the literature on OT planning and scheduling occurs on a daily basis. Long-term scheduling is convenient for surgeons and patients, however changes in waiting lists due to emergencies can cause difficulties. In this paper, we address this issue by considering planning and scheduling of the OTs over a two to four week period, whilst allowing for regular rescheduling to occur due to non-elective arrivals, cancellations, and changes in the waiting list.

---

*Author to whom all correspondence should be addressed (e-mail: b.spratt@qut.edu.au).

We present a robust Mixed Integer Programming (**MIP**) formulation of the integrated Master Surgical Schedule (**MSS**) Surgical Case Assignment (**SCA**) problem under a rolling horizon scheduling approach. Using a rolling horizon approach, we schedule for the upcoming planning horizon (a one to four week period), implement the first week of the schedule and shift the planning horizon forward. This approach also schedules operating room[1] (**OR**) time for patients, surgeons, and specialties. OR time blocks are reserved for non-elective arrivals. The sharing of ORs amongst both elective and non-elective patients is expected to improve schedule performance in terms of OR utilisation and patient waiting times (Wullink et al., 2007).

Whilst a case study of a large Australian public hospital motivates this approach, we do not incorporate Australian waiting list management policies, as these tend to perform poorly under a longer planning horizon. Instead, using metaheuristics, it is ensured that both their urgency category and the length of time they have spent on the waiting list prioritise otherwise identical patients. In doing so, we assign surgeons their category one patients in descending order of time-spent waiting, followed by their less urgent patients, also in descending order of time waiting.

The work presented in this paper is innovative in a number of ways. Firstly, the approach presented considers a significant level of detail to produce realistic and implementable schedules. The approach incorporates volatility in the surgical durations of both elective and non-elective patients, as well as the arrival of non-elective patients. Currently many hospitals use separate ORs for elective and non-elective patients. In this study, we also investigate the assignment of elective and non-elective OR time blocks. To account for cancellations and the changing waiting lists, a rolling horizon scheduling approach that limits deviations from an initial schedule is used. We use a rolling horizon approach to ensure hospital staff workload is predictable over a multi-week horizon, whilst incorporating that any additional surgery requests and cancellations in the surgical schedule. The literature contains very few instances of weekly reactive rescheduling of OTs. The metaheuristics used in this case are also quite innovative, particularly the development and subsequent use of new hyper metaheuristics on OT planning and scheduling problems.

The structure of the paper is as follows. In Section 2 we review the relevant literature, identify gaps, and discuss the novelty of the proposed approach. In Section 3 we present a description of the case study hospital and the problem addressed herein. We present the MIP formulation in Section 4. We list the assumptions made in Section 4.2. In Section 5, we present the solution techniques including details of the constructive heuristic and hyper metaheuristics implemented. Section 6 contains computational results. We discuss conclusions and perspectives in Section 7.

---

[1] An operating theatre (OT) is a set of operating rooms (ORs).

## 2. Literature Review

Throughout the literature, Operations Research is applied in a hospital setting to improve staff rostering (e.g. Rahimian et al. (2017)), perform capacity analysis (e.g. Burdett et al. (2017)), manage bed spaces (e.g. Landa et al. (2018)), and improve individual departments (e.g. Bai et al. (2018)). One such department, that has received much attention in recent years, is the OT department. Classifications of OT planning and scheduling problems are under three main levels: strategic, tactical, and operational. In this paper, we focus on the tactical and operational levels of OT planning and scheduling.

The main problem at the tactical level of OT planning and scheduling is the Master Surgical Scheduling Problem (**MSSP**). The MSSP is the problem of allocating surgeons, surgical teams, or surgical specialties OT time in the form of a MSS. This can be done under an open, block, or modified-block scheduling strategy. An open scheduling policy allows the scheduler to allocate surgeons to any suitable OR at any time. A block scheduling policy assigns entire blocks of time to surgeons or specialties. A modified block scheduling policy allows for a greater amount of flexibility in the allocation of time blocks.

Visintin et al. (2016) investigate how the flexibility of surgical teams, ORs and surgical units can improve the total number of patients scheduled using a MIP formulation of the MSSP. Allowing surgical teams to change with every new MSS and allowing a mix of both short and long stay cases in the same sessions maximised the number of patients scheduled. We utilise flexibility in our model (cf. Section 4) through a rolling horizon MSS that is updated to changing patient demands, whilst simultaneously making decisions at the operational level of OT planning and scheduling.

The two main problems at the operational level are the advanced scheduling and allocation scheduling of the ORs. Advanced scheduling involves determining which patients to treat each day (or time block) in the scheduling horizon. The Surgical Case Assignment Problem (**SCAP**) is another term for advanced scheduling. When solving the SCAP, the scheduling horizon is usually a single week. The second problem at the operational level, allocation scheduling, often occurs on a daily basis. Allocation scheduling, or the Surgical Case Sequencing Problem (**SCSP**), is the sequencing of patients within each day or block.

As each of the OT planning and scheduling problems are NP-hard, metaheuristics and hybrid metaheuristics are often used to produce good feasible solutions in reasonable amounts of computational time. Decomposition approaches are popular in the literature, however an integrated approach can lead to better solutions. Whilst the SCAP and SCSP are often addressed simultaneously (e.g. Molina-Pariente et al. (2015a), Moosavi and Ebrahimnejad (2018)), the MSSP and SCAP are mainly solved through decomposition approaches (e.g. Agnetis et al. (2014)). In this paper, we solve the integrated MSS SCA problem using a rolling horizon approach with consideration of the inherent uncertainty in the OT environment. A rolling horizon approach can improve schedule predictability for staff and patients, whilst updating to any changes in the waiting list, staff availability, and patient cancellations.

The inclusion of volatile elements in model formulations is another trend seen in the literature as it enables the production of schedules more robust to uncertainty. Van Riet and Demeulemeester (2015) review existing literature on the trade-offs involved in accounting for non-elective surgeries. The main objectives to consider are waiting time, utilisation, overtime, and cancellations. Ferrand et al. (2014) also review the literature on the ways to balance efficiency and responsiveness of the ORs in the face of non-elective demand.

To account for volatility, we consider lognormally distributed surgical durations and non-elective arrivals. This model maximises the number of elective surgeries performed, whilst reducing overtime and cancellations. To adapt to the changing waiting list we apply a rolling horizon scheduling approach. In the literature, daily rescheduling with a rolling horizon model produces less OR idle time and increases utilisation compared to a non-rolling horizon model (Luo et al., 2016). In their recent review of OT planning and scheduling, Zhu et al. (2018) indicate that there is a definite need for longer-term (several weeks) scheduling approaches (as opposed to daily approaches) that incorporate volatility. Our rolling horizon approach incorporates volatility in surgical durations, non-elective arrivals, and changes to the elective surgery waiting list on a longer term basis.

In our approach, we produce a robust MSS and SCA over periods of one to four weeks and reschedule every week, limiting the number of deviations from the original schedule. The MSS SCA is robust in that we reserve capacity for non-elective patients and consider lognormally distributed surgical durations such that the probability that a time block requires overtime is less than 5%. A four week MSS is in line with the case study hospital's current MSS update frequency. The case study hospital produces SCAs weekly. Thus, considering scheduling horizons of between one and four weeks is compliant with the hospital's current procedure.

Throughout the literature, there is still a need for more realistic models, with fewer simplifications, for application to hospital settings. There is also the need for solution methodologies that can handle long waiting lists and large surgical departments. We discuss hyper metaheuristics that produce good feasible solutions in reasonable amounts of computational time. Importantly, despite the large case study considered, standard desktop computers produce good solutions.

This paper contributes to the literature in a number of ways. Firstly, the model is realistic enough to account for a variety of staff and resource constraints and is implementable in the hospital environment to improve the current scheduling methodology. The model maximises surgical throughput whilst adhering to resource limitations. By considering volatile surgical durations and non-elective arrivals, we are not only able to produce a robust schedule that reduces cancellations and overtime, but using the rolling horizon scheduling approach, we are able to reschedule patients in the case of cancellations. In doing so, patient outcomes improve and overall satisfaction increases.

The MIP model presented utilises the inherent symmetry of the problem to reduce the number of constraints required. We successfully implement innovative hyper metaheuristics on an OT planning and scheduling model. These hyper metaheuristics produce better solutions than the comparative baseline metaheuristic. Finally, despite the size of the case study, we find good solutions in short amounts of computational time and the work is therefore accessible to hospital administrators.

## 3. Problem Description

In this paper, we address the problem of scheduling patients, surgeons, and specialties in ORs across a multiple weeks. The model uses a block scheduling policy, with either two half-day blocks or one full-day block. We simultaneously consider the MSS and SCA problems.

The problem pertains to a case study of a large Australian public hospital. Each year there are over 100,000 admissions to the hospital, with almost 70,000 admissions to the emergency department. There are around 15,000 elective surgeries each year, along with approximately 6,000 non-elective (emergency and urgent) surgeries. There is a long waiting list for elective surgeries currently composed of almost 2,900 patients.

Elective surgery requests are categorised by urgency and assigned a recommended waiting time of 30, 90, or 360 days. At present, the hospital is able to treat 100% of category one (the most urgent) patients within 30 days of being placed on the waiting list. This falls to around 90% for category two and three patients. We wish to improve the proportion of patients treated on time, by increasing OT throughput without increasing surgical overtime.

There is capacity for 825 bed spaces at the hospital, around 300 of these being surgical beds. The surgical department has 21 ORs. On weekdays, one OR is reserved for non-elective surgeries, although these tend to overflow into the other ORs. Administrative staff typically schedule elective surgeries on weekdays and never after eight o'clock at night. Four non-elective ORs are available on weekends. There are around 20 beds in both the surgical care unit and post anaesthesia care unit.

The case study hospital is a tertiary teaching hospital which results in procedures performed by staff of varying levels of experience. There is limited consultant availability, as many staff members prefer to spend the majority of their time working in private hospitals. The hospital can host visiting medical officers (VMOs) from other hospitals and can reserve OR time for these surgeries. The hospital outsources surgeries to surrounding hospitals if there is diminished theatre capacity or an increase in demand for a particular specialty. These intricacies make modelling the OT department quite a complicated process.

## 4. The Model

The model presented below relates to a case study of an Australian public hospital with a large surgical department. The model is detailed enough to produce implementable schedules in the hospital under study. The model is also generic enough that simple modifications to parameters and constraints will provide a model suitable for use in other hospitals.

The integrated MSSP and SCAP is formulated using MIP to allocate surgeries, surgeons and specialties to time blocks over a horizon comprised of several working weeks. In conjunction with the use of a rolling horizon scheduling approach, the inclusion of non-elective patients within the schedule is also considered.

Given the size of the case study, a commercial solver cannot solve this model (or any other formulation of the integrated MSS SCA) to optimality. The model presented here is NP-hard as both the MSS and SCA are reducible to bin-packing problems. Focus is placed on representing the problem, as opposed to reducing the number of constraints, as many of these would be eliminated by the pre-solve function of commercial solvers. For example, Constraints (17) and (16) could use index sets rather than binary parameters.

We include a wide range of real life constraints including those on surgeon availability and suitability. Surgeon availability data is provided in [dataset] Spratt and Kozan (2017). As this is a teaching hospital, it is important that only appropriate surgeons perform surgeries. The hospital prohibits surgeons from sharing time blocks. We consider OR suitability when allocating specialties to ORs. For example, we limit the movement of expensive equipment and allocate specialties accordingly.

The model reserves capacity for a certain number of non-elective patients of each specialty within the schedule. Hospital staff can adjust this to allow for a variety of risk attitudes. We assume that both elective and non-elective surgical durations are lognormally distributed (Spratt et al., 2018). Elective and non-elective patients cannot share OR blocks. A limited number of ORs can be opened on the weekend (based on surgeon availability) to accommodate non-elective surgeries. Elective surgeries must occur on weekdays. The model presented determines the optimum reservation of ORs for non-elective patients.

Due to the availability of data, we group surgical durations by specialty. The model does include individual duration parameters for each patient in the case that such data is available. In the case of specialties with highly variable durations, we split these into sub specialties based on patient urgency category.

Patients of the same specialty can share a time block as long as the $95^{th}$ percentile of their combined surgical durations is at most the length of the block. In the case where a single surgery's $95^{th}$ percentile exceeds the length of a full day block, the surgery can still take place in a full day block, as long as there are no other surgeries in that block. The lognormal parameters used herein are given in [dataset] Spratt and Kozan (2017). We refer to specialties by number rather than name to respect the confidentiality of hospital data.

Although we investigate disruptions by incorporating volatile surgical durations and non-elective arrivals, we apply a rolling horizon scheduling approach to reschedule in the case of disruptions. In doing so, we are able to update the schedule to account for cancellations and new arrivals on the waiting list. By producing an initial schedule over a four-week scheduling horizon, we are able to provide some insight into the future workload of hospital staff. This four-week scheduling horizon is in line with the current four-week Master Surgical Schedule used by the case study hospital. The rolling horizon scheduling approach allows for the flexibility required to create better schedules.

### 4.1. MIP Model Formulation

In this section, we present an MIP formulation for the integrated MSS SCA problem under a rolling horizon scheduling approach.

#### 4.1.1. Scalar Parameters

$\bar{H}$    the number of surgeons that practice at the hospital

| | |
|---|---|
| $\bar{P}$ | the number of patients in the waiting list at the start of the scheduling horizon |
| $\bar{S}$ | the number of surgical specialties |
| $\bar{R}$ | the number of ORs |
| $\bar{T}$ | the number of time periods in the scheduling horizon |
| $\bar{W}$ | the number of weeks in the scheduling horizon |
| $M_\psi$ | the maximum number of non-elective patients (of any specialty) that can be seen in a single OR in a full-day block |
| $M_p$ | the maximum number of elective patients (of any specialty) that can be seen in a single OR in a full-day block |
| $\xi$ | the maximum number of ORs that can be open at any one time during the weekend |

*4.1.2. Index Sets*

| | |
|---|---|
| $H$ | the set of surgeons that practice at the hospital. $H = \{1, \ldots, \bar{H}\}$ |
| $P$ | the set of patients on the waiting list at the start of the scheduling horizon. $P = \{1, \ldots, \bar{P}\}$ |
| $S$ | the set of surgical specialties. $S = \{1, \ldots, \bar{S}\}$ |
| $R$ | the set of ORs. $R = \{1, \ldots, \bar{R}\}$ |
| $T$ | the set of time periods in the scheduling horizon. $T = \{1, \ldots, \bar{T}\}$ |
| $W$ | the set of weeks in the scheduling horizon. $W = \{1, \ldots, \bar{W}\}$ |

*4.1.3. Indices*

| | |
|---|---|
| $h$ | index for surgeon in set $H$ |
| $p$ | index for patient in set $P$ |
| $r$ | index for OR in set $R$ |
| $s$ | index for surgeon in set $S$ |
| $t$ | index for time period in set $T$ |
| $w$ | index for week in set $W$ |

*4.1.4. Vector Parameters*

| | |
|---|---|
| $E_{ph}$ | 1 if patient $p$ can be treated by surgeon $h$, 0 otherwise, $\forall p \in P, h \in H$ |
| $F_{ht}$ | 1 if surgeon $h$ is available during time period $t$, 0 otherwise, $\forall h \in H, t \in T$ |
| $G_{hs}$ | 1 if surgeon $h$ is a member of specialty $s$, 0 otherwise, $\forall h \in H, s \in S$ |
| $I_{ps}$ | 1 if patient $p$ can be treated by specialty $s$, 0 otherwise, $\forall p \in P, s \in S$ |
| $R_{rs}$ | 1 if OR $r$ is equipped for surgeries by specialty $s$, 0 otherwise, $\forall r \in R, s \in S$ |
| $D_{t\tau}$ | 1 if time periods $t$ and $\tau$ do not overlap, 0 otherwise, $\forall t, \tau \in T$ |
| $B_t$ | 1 if time period $t$ is a full day block, 0 otherwise, $\forall t \in T$ |
| $V_t$ | 1 if time period $t$ is on the weekend, 0 otherwise, $\forall t \in T$ |
| $U_{tw}$ | 1 if time period $t$ is in week $w$, 0 otherwise, $\forall t \in T, w \in W$ |
| $\kappa_s^+$ | the number of specialty $s$ elective patients treatable in a full day block, $\forall s \in S$ |
| $\kappa_s^-$ | the number of specialty $s$ elective patients treatable in a half day block, $\forall s \in S$ |
| $\hat{\kappa}_s^+$ | the number of specialty $s$ non-elective patients treatable in a full day block, $\forall s \in S$ |

$\hat{\kappa}_s^-$  the number of specialty $s$ non-elective patients treatable in a half day block, $\forall s \in S$

*4.1.5. Decision Variables*

$X_{srt}$   1 if specialty $s$ is assigned to OR $r$, time period $t$, 0 otherwise, $\forall s \in S, r \in R, t \in T$

$Y_{hrt}$   1 if surgeon $h$ is assigned to OR $r$, time period $t$, 0 otherwise, $\forall h \in H, r \in R, t \in T$

$Z_{prt}$   1 if patient $p$ is treated in OR $r$, time period $t$, 0 otherwise, $\forall p \in P, r \in R, t \in T$

$\Psi_{srt}$   the number of non-elective specialty $s$ patients for OR $r$, period $t$, $\forall s \in S, r \in R, t \in T$

*4.1.6. Objective Function*

The objective of the rolling horizon scheduling model is to maximise the number of elective surgeries performed during the scheduling horizon. This will have the effect of reducing the waiting list, decreasing wait time, and improving patient outcomes.

$$\text{Maximise} \sum_{p \in P} \sum_{r \in R} \sum_{t \in T} Z_{prt} \tag{1}$$

*4.1.7. Constraints*

Constraint (2) ensures that we treat each patient at most once during the scheduling horizon.

$$\sum_{r \in R} \sum_{t \in T} Z_{prt} \leq 1, \quad \forall p \in P \tag{2}$$

We assign no more than one specialty (3) and no more than one surgeon (4) to an OR during a time block or overlapping time block. By enforcing these constraints for $t < \tau$, the symmetry of the constraint is utilised.

$$\sum_{s \in S} X_{srt} + X_{sr\tau} \leq 1 + D_{t\tau}, \quad \forall r \in R, t, \tau \in T, t < \tau \tag{3}$$

$$\sum_{h \in H} Y_{hrt} + Y_{hr\tau} \leq 1 + D_{t\tau}, \quad \forall r \in R, t, \tau \in T, t < \tau \tag{4}$$

If the model assigns a surgeon to a time period, then it must assign the correct specialty to that time period, in that OR.

$$Y_{hrt} \leq \sum_{s \in S} G_{hs} X_{srt}, \quad \forall h \in H, r \in R, t \in T \tag{5}$$

A specialty can only be assigned to an OR if the OR is equipped for that specialty.

$$X_{srt} \leq R_{rs}, \quad \forall s \in S, r \in R, t \in T \tag{6}$$

If the model assigns a patient to a time period, then it must assign the correct surgeon (constraint (7)) and specialty (constraint (8)) to that time period, in that OR.

$$Z_{prt} \leq \sum_{h \in H} E_{ph} Y_{hrt}, \quad \forall p \in P, r \in R, t \in T \tag{7}$$

$$Z_{prt} \leq \sum_{s \in S} I_{ps} X_{srt}, \quad \forall p \in P, r \in R, t \in T \tag{8}$$

If we reserve an OR block for non-elective capacity, then we must assign a surgeon to that time period. The surgeon assigned must be a member of the appropriate specialty.

$$\Psi_{srt} \leq M_\psi \sum_{h \in H} Y_{hrt} G_{hs}, \quad \forall s \in S, r \in R, t \in T \tag{9}$$

The model can only assign available surgeons.

$$\sum_{r \in R} Y_{hrt} \leq F_{ht}, \quad \forall h \in H, t \in T \tag{10}$$

Based on the 95th percentile of the sum of lognormal surgical durations fitting within a block, we restrict the number of elective (constraint (12)) and non-elective (constraint (13)) patients assigned to each OR block.

$$\sum_{p \in P} Z_{prt}\, I_{ps} \leq \left(\kappa_s^+ B_t + \kappa_s^- (1 - B_t)\right) X_{srt}, \qquad \forall s \in S, r \in R, t \in T \tag{11}$$

$$\Psi_{srt} \leq \left(\hat{\kappa}_s^+ B_t + \hat{\kappa}_s^- (1 - B_t)\right) X_{srt}, \qquad \forall s \in S, r \in R, t \in T \tag{12}$$

An OR block can be reserved for either elective patients or non-elective patients, but cannot be shared.

$$\frac{1}{M_\psi} \sum_{s \in S} \Psi_{srt} \leq 1 - Z_{prt}, \qquad \forall r \in R, t \in T, p \in P \tag{13}$$

Surgeons cannot perform elective surgery on the weekend.

$$\sum_{p \in P} \sum_{r \in R} Z_{prt} \leq M_p (1 - V_t), \qquad \forall t \in T \tag{14}$$

Each week, we reserve capacity for a certain number of non-elective patients of each specialty.

$$\sum_{r \in R} \sum_{t \in T} \Psi_{srt} U_{tw} \geq \psi_s, \qquad \forall w \in W, s \in S \tag{15}$$

To ensure that we do not rely on weekends too heavily for non-elective surgeries, we limit the number of ORs open at any one time throughout the weekend.

$$\sum_{r \in R} \sum_{s \in S} \left( X_{srt} + \sum_{\tau \in T} D_{t\tau} X_{sr\tau} \right) \leq \xi + \overline{R}(1 - V_t), \qquad \forall t \in T \tag{16}$$

4.2. The Rolling Horizon Approach

When considering OT planning problems, in particular the MSSP and SCAP, most authors produce schedules for a single week with no regard for future demand. Rolling horizon decision making is the process of making decisions based on predictions of the future, and adjusting to deviations as they occur.

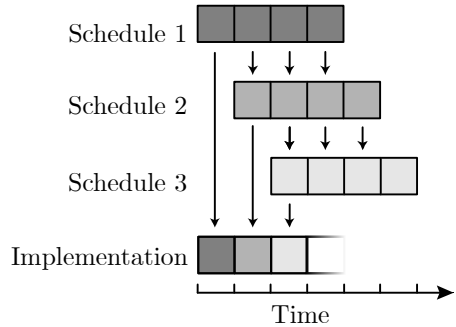

**Figure 1:** A rolling horizon approach to scheduling.

Using a rolling horizon approach, we produce a schedule for the length of the planning horizon. The hospital implements only the first week of the schedule and updates the data with any new arrivals and cancellations. We then produce a new schedule for the upcoming planning horizon. We repeat this process until the hospital has implemented schedules for the entire planning period. Figure 1 shows this process.

At the end of each week, hospital staff must update the data before rescheduling. The parameters that require regular updates are the waiting list related parameters ($\bar{P}, P, E_{ph}, I_{ps}$). The waiting list related parameters will change whenever a patient cancels or if the hospital receives a new elective surgery request. Other parameters may change if, for example, a surgeon updates their availability or an OR is under maintenance.

When applying a rolling horizon scheduling approach to the model presented above, additional constraints are required, limiting the number of deviations from the original schedule. Here, the only schedule deviation we consider is the change in a patient's surgery time period. We include the previous patient schedule as $Z^*_{prt}$, the new decision variable, $\nu_p$, and scalar parameter $\rho$.

$Z_{prt}$    1 if patient $p$ is treated in OR $r$, time period $t$, 0 otherwise, $\forall p \in P, r \in R, t \in T$

$\nu_p$    1 if patient $p$ was scheduled previously, in a different time period, $\forall p \in P$

$\rho$    maximum proportion of patients deviations in the first $\overline{W} - 1$ weeks of the new schedule

Constraint (17) ensures that if a patient was scheduled for time period $t$ in the previous schedule, but is not scheduled for time period $t$ in the current schedule, then $\nu_p = 1$.

$$\sum_{r \in R} Z^*_{prt} - \sum_{r \in R} Z_{prt} \leq \nu_p, \forall p \in P, t \in T \setminus \{\overline{W}\} \quad (17)$$

Constraint (18) ensures that the correct proportion of patients retain their previously scheduled surgical date.

$$\sum_{p \in P} \nu_p \leq \rho \sum_{p \in P} \sum_{r \in R} \sum_{t \in T \setminus \{\overline{W}\}} Z^*_{prt} \quad (18)$$

*4.3. Model Assumptions*

In this subsection we discuss the assumptions made when implementing the MSS SCA model on the case study. These assumptions are based on a case study of a large Australian public hospital. We assume that historical data is an accurate representation of the number of ORs that can be open, surgeon availability and specialty, surgical duration, and non-elective arrivals.

A modified block scheduling policy is used as per current hospital policy. This means that patients of multiple specialties cannot be treated in the same block. Any ORs that are not reserved under the MSS can be reallocated or used for non-elective surgeries. No ORs are dedicated to non-elective surgery, but OR blocks may be dedicated to either non-elective or elective surgery. Non-elective patients are prioritized over elective patients. As such, it is assumed that non-elective patients are treated in the first available time block.

Working days are ten hours. We assume that there is no flexibility in the working day, however if overtime is required a surgeon will perform the overtime. If a surgeon is not available for a particular shift, that surgeon will not be scheduled. The surgeon named on the waiting list is the surgeon that must be present during the surgery. If a surgeon is named on the waiting list, then that surgeon may not supervise additional surgeries in different theatres at the time of surgery.

The suitability of ORs must be respected. If a particular specialty does not require specialized equipment, then that surgery is listed as a 'general' surgery. We also assume that sufficient capacity exists in downstream wards, ensuring that bottlenecks are highly unlikely.

## 5. Solution Approaches

The model formulated in Section 4 is NP-hard. As such we require the use of metaheuristics to produce solutions in reasonable amounts of computational time. In the search for efficient solution strategies, heuristics and metaheuristics are often hybridised with other solution techniques, including mathematical programming.

Hyper metaheuristics are a class of solution techniques that use heuristic methodology when either selecting or generating heuristics (Burke et al., 2013). Unlike other heuristics and metaheuristics, the search space of hyper metaheuristics is a set of heuristics, rather than the solution space. For example, neighbourhood swaps may be selected randomly by a given metaheuristic. A hyper metaheuristic version may use a heuristic or metaheuristic for selecting the neighbourhood swap type to implement at each step of the algorithm.

Throughout the literature, it is common to use hyper metaheuristics to solve combinatorial optimisation problems including timetabling and scheduling. Hyper metaheuristics are uncommon in hospital-planning literature. Although a number of problems in health care planning use hyper metaheuristics, no authors have applied hyper metaheuristics to the integrated MSS SCA problem. The development and use of new hyper metaheuristics on the integrated MSS SCA problem is an innovative contribution to the literature. These hyper metaheuristics provide promising results. There are no constructive heuristics for the MSS SCA problem in the literature. An efficient constructive heuristic could improve hospital performance by providing good solutions quickly. This approach can also provide scenario analysis and plan potential capacity changes.

In this paper, we develop new hyper metaheuristics for use on the integrated MSS SCA. We implement Hyper Simulated Annealing (**SA**), which selects neighbourhood type through use of a Tabu Search (**TS**) strategy and accepts solutions using a SA approach. We implement Hyper SA-TS using the Hyper SA algorithm along with a Tabu list of moves. By implementing each of these algorithms, we are able to utilise the problem structure whilst obtaining good solutions in reasonable amounts of computational time. We use SA to compare solution quality and computational effort with the proposed Hyper SA and Hyper SA-TS metaheuristics.

When running the metaheuristics, we consider scheduling horizon lengths of between one and four weeks. Scheduling horizons of between one and four weeks are compliant with the hospital's current procedure, and can have a large impact on scheduling objectives (Molina-Pariente et al., 2015b).

We run the metaheuristic, providing a MSS and SCA for the entire scheduling horizon. Stepping forward a week, we remove patients treated in the first week from the waiting list, and update the list with any new arrivals. The schedule calculated in the last iteration is the initial schedule for the metaheuristic in the next iteration.

### 5.1. Baseline Non-elective Schedule

Given the complexity of ensuring sufficient OR capacity is reserved for non-elective surgeries, in this subsection we present a reduced version of the model provided in Section 4. This reduced model is used to allocate non-elective OR reservations and is implemented before performing the elective scheduling. We do this in CPLEX. We make a simplifying assumption when creating the baseline non-elective schedule. We assume that staff treat non-elective patients in the first reserved block and prioritise non-elective patients above elective patients.

To maintain alignment with the original objective of maximising the number of elective patients treated, we create the baseline non-elective schedule using the objective shown in (19). This objective is to minimise the amount of OR time reserved for non-elective surgeries.

$$\text{Minimise } \sum_{s \in S} \sum_{r \in R} \sum_{t \in T} X_{srt} \tag{19}$$

In solving this problem, we minimise (19) subject to constraints (2) to (10), (12), (15) and (16). Fortunately, CPLEX is able to solve this problem in a reasonable amount of time. We use this initial solution when solving for the elective portion of the schedule using metaheuristics. Although tackling the non-elective and elective portions of the schedule separately may reduce solution quality, any reduction in solution quality is negligible compared to the computational effort required to create a full schedule. While it is possible to vary this baseline non-elective schedule, it is unnecessary to do so. Each week we use the same baseline schedule, and as such, the time taken to compute the initial solution is unimportant.

## 5.2. Constructive Heuristic

Each of the metaheuristics discussed previously are adjustable to suit a myriad of problems through simple changes to the solution neighbourhoods. As such, the metaheuristics cannot fully utilise known problem structure. In this subsection, we develop a constructive heuristic for the initial solution of the MSS SCA problem.

Since each patient is preassigned a surgeon and specialty, and surgical durations parameters are the same within each specialty (or sub specialty), we group the waiting list into unique surgeon-specialty combinations. Each surgeon-specialty set has a number of characteristics:

- the number of patients on the waiting list, if the surgeon is available for a full day block on a weekday,
- the maximum number of patients (based on lognormal surgical durations, derived from historical data) that can be treated in either a full day or half-day block (as appropriate), and
- the actual maximum number of patients in a single block.

The actual maximum number of patients is the minimum of:

- the number of patients that could be treated by that surgeon-specialty according to surgical durations
- the number of patients remaining on the waiting list for that surgeon-specialty

If a surgeon-specialty combination has, at most, enough patients waiting to fill two half-day blocks, we do not assign the surgeon a full-day time block and adjust their availability to reflect this. This is as, for each specialty, the number of patients that fit into a full day block is at least two times the number of patients that fit into a half day block. By limiting the allocation of full day time blocks in this way, we maintain schedule flexibility and do not allocate full day blocks unless they can be fully utilised.

When selecting which surgeon-set to schedule first, we sort the actual maximum number of patients in a block in descending order; this is a greedy approach similar to real-life. Ties are broken by considering the average maximum number of patients (according to lognormal surgical durations) in ascending order, surgeons unavailable for full day blocks before available surgeons, and the total number of patients waiting for treatment by each surgeon-specialty set. In sorting in this order, we fill the schedule using a greedy approach. The consideration of regular surgeon-specialty capacity ensures the model prioritises surgeon-specialty sets that can fill their allocated time blocks. This enables new patients on the waiting list to complete partially full surgical lists in the following weeks. By breaking further ties through consideration of full day availability, the model fills half-day blocks that could hold as many surgeries as full day blocks first. The final tiebreak, using total number of patients waiting, aims to reduce the longer surgical waiting lists by essentially assigning a higher priority to those surgeon-specialty sets.

After selecting the surgeon-specialty set, we schedule the surgeon-specialty set in every time block available. When considering which OR to schedule the surgeon-specialty set **a regret-based approach** is used. For each OR, we calculate the actual maximum number of patients (based on the remaining waiting list) for each feasible surgeon-specialty set.

The regret for not scheduling the surgeon-specialty set in the OR is the difference between:

- the highest number of patients that could to be scheduled under the selected surgeon-specialty set, and
- the second highest number of patients that could be scheduled under all feasible surgeon-specialty sets

The OR selected is the one with the largest regret. The model breaks ties by considering the third best actual maximum under all feasible specialty sets.

Not only is this constructive heuristic computationally efficient, but the use of a greedy heuristic combined with a regret based OR selection produces good feasible solutions to the MSS SCA problem. Figure 2 shows the flow diagram.

*5.3. Hyper SA*

Like hybridised metaheuristics, hyper metaheuristics often produce better solutions than the original metaheuristics in shorter times. Whilst hybridised metaheuristics feature prominently in the literature, hyper metaheuristics are far less common. Hyper metaheuristics select the local search heuristics to use at each step based on a metaheuristic approach whilst using another metaheuristic approach determines whether to accept the solution updates.

A metaheuristic's choice of solution neighbourhood (heuristic) has a big impact on solution quality. To exploit this fact, we implement a hyper metaheuristic that uses past neighbourhood performance to determine the next neighbourhood to select. We use TS to select the heuristic to use at each iteration within the SA metaheuristic.

This hyper metaheuristic uses TS in conjunction with heuristic rankings to select the heuristic for each set of iterations. A SA approach accepts solutions such that improving solutions are always accepted and the algorithm accepts worsening solutions according to the change in solution and the solution temperature. In addition to the standard temperature decrease, the algorithm performs temperature increases throughout the run. To account for fewer improving moves towards the end of a run, the algorithm performs iterations of the same heuristic type in blocks when calculating heuristic rank.

In summary, the algorithm selects the neighbourhood swap type with the highest rank out of those that are not tabu. Using this neighbourhood swap type, the algorithm generates a random neighbouring solution. If this solution is feasible and acceptable according to the acceptance probability function, we update the solution and reduce the algorithm temperature. Otherwise, the solution remains unchanged and the algorithm temperature is increased. After a predefined number of iterations, we update the tabu list.

The Appendix contains details of the algorithm. Figure 3 shows the Hyper SA flow diagram. Figure 4 shows the Tabu Update subroutine.

### 5.4. Hyper SA-TS

After initial investigation showed that the use of TS significantly improved performance in certain problem instances, we include a list of tabu moves when implementing Hyper SA-TS. A list of the most recent theatre-week-day combinations is stored. The algorithm cannot perform a neighbour swap on a tabu theatre-week-day combination. If a move worsens the solution, we add theatre-week-day combination to the tabu list. If the tabu list exceeds a given length, the algorithm removes the oldest tabu move. Every N iterations the maximum allowable tabu length is decreased. This gives the metaheuristic more 'freedom' in block selection towards the end of the run. The tabu theatre-week-day combinations are entirely separate from the list of tabu neighbourhood types, denoted Tabu_twd and Tabu_n respectively. More sophisticated implementations are possible; however this methodology performs well for this problem. Figure 5 shows the flow diagram (see Appendix).

The use of a tabu list of neighbourhood swap types enables the algorithm to select the neighbourhood swap type of the highest rank, of those that are non-tabu. We add a neighbourhood swap type to the list of tabu moves if after a predefined number of iterations it has not improved the solution. If, after a predefined number of iterations, the solution has improved, the rank of the solution neighbourhood type is increased. Additionally, the use of an acceptance probability function ensures that worsening moves are occasionally accepted, enabling the algorithm to escape local optima.

After a good neighbourhood swap, the temperature decreases, reducing the probability that we accept worsening swaps in the future. Conversely, after a worsening neighbourhood swap, solution temperature increases, increasing the probability that we accept worsening swaps. Through this methodology, if it is likely that the algorithm is stuck in local optima (and therefore unable to find improving solutions through simple neighbourhood swaps) it is more likely that the algorithm will accept a worsening solution and move away from the local optima.

## 6. Results

In this section, we present the results of computational experiments using SA, Hyper SA, and Hyper SA-TS on a case study of a large Australian public hospital. The surgical department under study is quite large, with 21 ORs, over 100 surgeons, and 27 specialties (including subspecialties required under the formulation presented in Section 4). At present, the hospital has around 2,900 patients on the elective surgery waiting list. The hospital treats approximately 15,000 elective surgery patients per year and a further 6,000 non-elective surgical patients.

We base the computational instances in this paper on the case study hospital. In these computational experiments, we consider 108 surgeons, 21 ORs, 27 specialties, 21 time blocks, and an initial waiting list of 2871 patients. The baseline non-elective schedule reserves capacity for 113 non-elective patients each week. Details regarding the surgeon availability, surgical durations, initial waiting list and average number of surgical requests can be found in [dataset] Spratt and Kozan (2017).

We coded the metaheuristics in MATLAB® and average performance was determined after 100 runs on the university's High Performance Computing (**HPC**) facility. The solution variance had converged sufficiently after 100 runs of each metaheuristic. To demonstrate the ease of implementation in a hospital environment, average computational time was found using MATLAB® on an Intel® Core™ i7-370 CPU @ 3.40GHz with 16 GB of RAM.

The maximum total number of iterations was limited to 16,000 in each case as computational time is a valuable resource in the hospital. Under a one-week scheduling horizon, 16,000 iterations were sufficient for metaheuristic convergence. We found the other parameters through parameter tuning performed on the university's HPC facility. In this section, the time shown is the average time required to create the schedules for the entire planning period (scheduling horizon lengths vary) each week for a total of six weeks of data available. We do not include the time required to calculate the baseline schedule (see Section 5.1) as we only generate the baseline schedule and only ever update it if the predicted non-elective surgery demand changes.

### 6.1. Varying Scheduling Horizon

Whilst in a static environment, full scheduling is superior to partial scheduling (Sabuncuoglu and Bayiz, 2000), in a dynamic and volatile environment there is merit to both approaches. The effect of horizon length on approach performance is analysed for the rolling horizon scheduling approach on this particular instance of the MSS SCA problem.

In this subsection, we consider the effect of varying scheduling horizon length on the mean total number of patients scheduled. We performed parameter tuning for each of the metaheuristics under each scheduling horizon length. The maximum total iterations were limited to 16,000 to keep computational time low whilst ensuring the metaheuristics converge.

Table 1 displays results from these computational experiments. Table 1 includes the mean total number of patients scheduled, the variance, worst and best objectives, and the time taken (in seconds) to schedule all six weeks of historical data, whilst implementing a rolling horizon scheduling approach with scheduling horizon lengths between one and four weeks. For example, under a four week scheduling horizon we produce a four week schedule, update waiting lists assuming the first week went to plan, and reschedule for another four week scheduling horizon until all six weeks of historical data are scheduled.

The results displayed in Table 1 indicate that Hyper SA is the best performing metaheuristic under scheduling horizons of one to four weeks. Hyper SA appears to outperform Hyper SA-TS under scheduling horizons of two ($p<<0.01$), three ($p<0.01$), and four ($p<0.05$) weeks.

**Table 1:** Computational Results

| Horizon Length (weeks) | Method | Mean | Variance | Worst | Best | Time (seconds) |
|---|---|---|---|---|---|---|
| 1 | SA | 1359.34 | 25.78 | 1346 | 1371 | 96.53 |
|   | **Hyper SA** | **1361.70** | **27.02** | **1348** | **1372** | **103.49** |
|   | Hyper SA-TS | 1361.63 | 23.57 | 1349 | 1373 | 113.29 |
|   | Constructive | 1284.00 | - | - | - | 1.49 |
| 2 | SA | 1354.09 | 42.24 | 1338 | 1368 | 99.44 |
|   | **Hyper SA** | **1365.85** | **21.16** | **1356** | **1376** | **129.71** |
|   | Hyper SA-TS | 1363.80 | 18.57 | 1352 | 1375 | 135.61 |
|   | Constructive | 1272.00 | - | - | - | 2.68 |
| 3 | SA | 1341.87 | 64.68 | 1320 | 1359 | 104.97 |
|   | **Hyper SA** | **1364.14** | **19.39** | **1353** | **1374** | **145.34** |
|   | Hyper SA-TS | 1362.52 | 22.27 | 1349 | 1374 | 153.17 |
|   | Constructive | 1205.00 | - | - | - | 3.30 |
| 4 | SA | 1310.18 | 92.03 | 1290 | 1330 | 109.79 |
|   | **Hyper SA** | **1346.97** | **36.39** | **1329** | **1360** | **153.04** |
|   | Hyper SA-TS | 1345.19 | 25.43 | 1335 | 1358 | 164.73 |
|   | Constructive | 1281.00 | - | - | - | 5.04 |

The actual number of elective patients treated across the six weeks was 1228, thus using the methodology presented in this paper we are able to clear a significant amount of backlog. Using the proposed metaheuristic approaches, we schedule approximately 24 additional patients each week whilst reserving more than enough capacity for non-elective patients. We do this while also ensuring number of ORs in use is in-line with historical averages. If the hospital were to consider a higher risk attitude toward overtime, it would be possible to schedule even more patients.

The proposed constructive heuristic outperforms historical hospital schedules under scheduling horizons of one, two, and four weeks. The constructive heuristic is simple to implement and runs quickly, providing hospital management with the opportunity to produce an initial schedule on a weekly basis, and perform modifications based on expert knowledge. Whilst the metaheuristic approaches outperform the proposed constructive heuristic, it is still a valuable planning tool for administrative staff members as it appears to mimic current scheduling strategy.

Two-sample t-tests indicate that Hyper SA (p<0.01) and Hyper SA-TS (p<0.05) performs best under a two-week scheduling horizon. SA performs best under a scheduling horizon of one week (p<<0.01). This may be due to the increase in solution space associated with a longer scheduling horizon, without any increase in maximum total iterations. To further investigate this, we perform parameter tuning for each of the metaheuristics with maximum total iterations now 16,000 times the length of the scheduling horizon (in weeks). Table 2 displays the new computational results.

When we increase the total number of iterations proportionally to scheduling horizon length, Hyper SA remains the top performing metaheuristic under scheduling horizons of one, two, and three weeks, whilst Hyper SA-TS is the top performing metaheuristic using a four-week scheduling horizon approach. Under scheduling horizons of three and four weeks, two-sample t-tests indicate that there may be no significant difference between Hyper SA and Hyper SATS (p>0.05).

Table 2 shows that the metaheuristics display a significant increase in number of patients scheduled when we increase the number of iterations proportionally to scheduling horizon length. Overlapping 95% confidence intervals and two-sample t-tests (p>0.05) show that under the increased total iterations Hyper SA and Hyper SA-TS may not perform significantly better under a scheduling horizon of three weeks compared to a two-week scheduling horizon. Under the increased total iterations, SA performs significantly better with a scheduling horizon of two weeks (p<0.01).

**Table 2:** Computational Results – Effect of Increasing Total Iterations

| Horizon Length (weeks) | Method | Mean | Variance | Worst | Best | Time (seconds) |
|---|---|---|---|---|---|---|
| 1 | SA | 1359.34 | 25.78 | 1346 | 1371 | 96.53 |
|   | **Hyper SA** | **1361.70** | **27.02** | **1348** | **1372** | **103.49** |
|   | Hyper SA-TS | 1361.63 | 23.57 | 1349 | 1373 | 113.29 |
| 2 | SA | 1365.45 | 26.51 | 1351 | 1376 | 190.50 |
|   | **Hyper SA** | **1369.34** | **16.69** | **1359** | **1379** | **232.17** |
|   | Hyper SA-TS | 1368.23 | 16.91 | 1360 | 1378 | 243.58 |
| 3 | SA | 1363.35 | 24.92 | 1344 | 1375 | 298.75 |
|   | **Hyper SA** | **1369.55** | **14.90** | **1360** | **1379** | **363.24** |
|   | Hyper SA-TS | 1368.87 | 13.29 | 1361 | 1379 | 393.82 |
| 4 | SA | 1347.03 | 37.36 | 1323 | 1358 | 411.41 |
|   | Hyper SA | 1359.22 | 22.92 | 1345 | 1369 | 487.58 |
|   | **Hyper SA-TS** | **1359.39** | **19.82** | **1346** | **1367** | **515.55** |

The increased throughput is due to an increase in OT utilisation and the choice of surgical specialty. The algorithm gives preference to surgical specialties with shorter surgeries. This is in line with the shortest first heuristic rules which tend to perform well on scheduling problems. For example, the algorithm allocates many time blocks to ophthalmology, as they are able to treat up to eight patients in a full-day block. After this initial clearing of backlog, we expect lower weekly throughput as bottlenecks become more apparent.

There is a significant reduction in overtime by considering the 95[th] percentile of surgical durations. Using the 95[th] percentile of lognormal surgical durations, we do not see any overtime caused by elective surgeries. Overtime may occur once throughout the week during the time reserved for liver transplants. This is as the 95$^{th}$ percentile of the duration of a single liver transplant is approximately 11.82 hours. Thus, the total overtime seen in the new schedule is only 1.82 hours. This is only 0.19% of the 950.54 hours scheduled.

In comparison, the historical schedule contains a total of 54.22 hours of surgery outside of 8am to 6pm. This is approximately 7.31% of the 741.80 surgical hours used throughout the week. This difference in overtime may be due to the difference in the reservation of capacity for non-elective surgeries. At present, the hospital dedicates OR8 and OR10 to non-elective surgeries. The majority of overtime occurs in these ORs. In the model presented in this paper, we do not reserve specific ORs for non-elective surgeries and instead reserve a variety of ORs throughout the week. In reality, surgeries may be too urgent to hold until the next morning and will occur as soon as possible. This would increase the overtime observed when implementing the new schedule.

## 7. Conclusion

In this paper we presented a MIP formulation of the combined MSSP SCAP. The model includes volatile surgical durations, non-elective arrivals and a rolling horizon scheduling approach. In this model, we were able to address the presence of non-elective patients by reserving OT capacity in blocks. In doing so, we also allowed non-elective surgeries to occur on weekends. The surgical durations of both elective and non-elective patients are lognormally distributed, and obtained from historical data. As such, constraints ensured that the 95$^{th}$ percentile of surgical durations does not exceed the length of the time block. In allocating time for surgeries, we ensure that elective and non-elective patients do not share time blocks.

We consider a number of solution techniques including metaheuristics and hyper metaheuristics. In particular we find that Hyper SA and Hyper SA-TS perform well compared to historical averages. We implemented a constructive heuristic designed specifically for the MSS SCA problem. The constructive heuristic outperforms historical schedules and requires only a fraction of the time.

One of the model limitations lies in the allocation of non-elective capacity. It is assumed that hospital administrators will shift elective patients to non-elective blocks if highly urgent cases arrive, allowing non-elective surgeries to be performed in the newly empty elective block. To address this limitation, future work includes producing and validating a number of simple heuristics for the reactive rescheduling of the ORs in the case of schedule disturbances.

The approach we presented in this paper is implementable in hospitals worldwide. By aligning the solution techniques with real-time information updates, it is possible to produce a decision support tool for administrative staff. The hyper metaheuristics solve the real-life complex problem quickly and can provide updated feasible schedules to hospital staff in a timely manner.

## Acknowledgements

This research was funded by the Australian Research Council (ARC) Linkage Grant LP 140100394. Computational resources and services used in this work were provided by the HPC and Research Support Group, Queensland University of Technology, Brisbane, Australia.

**References**

AGNETIS, A., COPPI, A., CORSINI, M., DELLINO, G., MELONI, C. & PRANZO, M. 2014. A decomposition approach for the combined master surgical schedule and surgical case assignment problems. *Health Care Management Science,* 17**,** 49-59.

BAI, J., FÜGENER, A., SCHOENFELDER, J. & BRUNNER, J. O. 2018. Operations research in intensive care unit management: a literature review. *Health Care Management Science,* 21**,** 1-24.

BURDETT, R. L., KOZAN, E., SINNOTT, M., COOK, D. & TIAN, Y.-C. 2017. A mixed integer linear programing approach to perform hospital capacity assessments. *Expert Systems with Applications,* 77**,** 170-188.

BURKE, E. K., GENDREAU, M., HYDE, M., KENDALL, G., OCHOA, G., ÖZCAN, E. & QU, R. 2013. Hyper-heuristics: A survey of the state of the art. *Journal of the Operational Research Society,* 64**,** 1695-1724.

FERRAND, Y. B., MAGAZINE, M. J. & RAO, U. S. 2014. Managing operating room efficiency and responsiveness for emergency and elective surgeries - A literature survey. *IIE Transactions on Healthcare Systems Engineering,* 4 (1)**,** 49-64.

LANDA, P., SONNESSA, M., TÀNFANI, E. & TESTI, A. 2018. Multiobjective bed management considering emergency and elective patient flows. *International Transactions in Operational Research,* 25**,** 91-110.

LUO, L., LUO, Y., YOU, Y., CHENG, Y., SHI, Y. & GONG, R. 2016. A MIP Model for Rolling Horizon Surgery Scheduling. *Journal of medical systems,* 40**,** 1-7.

MOLINA-PARIENTE, J. M., FERNANDEZ-VIAGAS, V. & FRAMINAN, J. M. 2015a. Integrated operating room planning and scheduling problem with assistant surgeon dependent surgery durations. *Computers & Industrial Engineering,* 82**,** 8-20.

MOLINA-PARIENTE, J. M., HANS, E. W., FRAMINAN, J. M. & GOMEZ-CIA, T. 2015b. New heuristics for planning operating rooms. *Computers & Industrial Engineering,* 90**,** 429-443.

MOOSAVI, A. & EBRAHIMNEJAD, S. 2018. Scheduling of elective patients considering upstream and downstream units and emergency demand using robust optimization. *Computers & Industrial Engineering,* 120**,** 216-233.

RAHIMIAN, E., AKARTUNALı, K. & LEVINE, J. 2017. A hybrid Integer Programming and Variable Neighbourhood Search algorithm to solve Nurse Rostering Problems. *European Journal of Operational Research,* 258**,** 411-423.

RAIS, A. & VIANA, A. 2011. Operations Research in Healthcare: a survey. *International Transactions in Operational Research,* 18**,** 1-31.

SABUNCUOGLU, I. & BAYIZ, M. 2000. Analysis of reactive scheduling problems in a job shop environment. *European Journal of Operational Research,* 126**,** 567-586.

SAMUDRA, M., VAN RIET, C., DEMEULEMEESTER, E., CARDOEN, B., VANSTEENKISTE, N. & RADEMAKERS, F. E. 2016. Scheduling operating rooms: achievements, challenges and pitfalls. *Journal of Scheduling,* 19 (5).

SPRATT, B. & KOZAN, E. 2017. Australian public hospital surgical case study. Mendeley Data.

SPRATT, B., KOZAN, E. & SINNOTT, M. 2018. Analysis of uncertainty in the surgical department: durations, requests, and cancellations. *Aust Health Rev.*

VAN RIET, C. & DEMEULEMEESTER, E. 2015. Trade-offs in operating room planning for electives and emergencies: A review. *Operations Research for Health Care,* 7**,** 52-69.

VISINTIN, F., CAPPANERA, P. & BANDITORI, C. 2016. Evaluating the impact of flexible practices on the master surgical scheduling process: an empirical analysis. *Flex Serv Manuf J,* 28**,** 182-205.

WULLINK, G., VAN HOUDENHOVEN, M., HANS, E. W., VAN OOSTRUM, J. M., VAN DER LANS, M. & KAZEMIER, G. 2007. Closing emergency operating roomsimproves efficiency. *Journal of Medical Systems,* 31**,** 543-546.

ZHU, S., FAN, W., YANG, S., PEI, J. & PARDALOS, P. M. 2018. Operating room planning and surgical case scheduling: a review of literature. *Journal of Combinatorial Optimization.*

**Appendix**

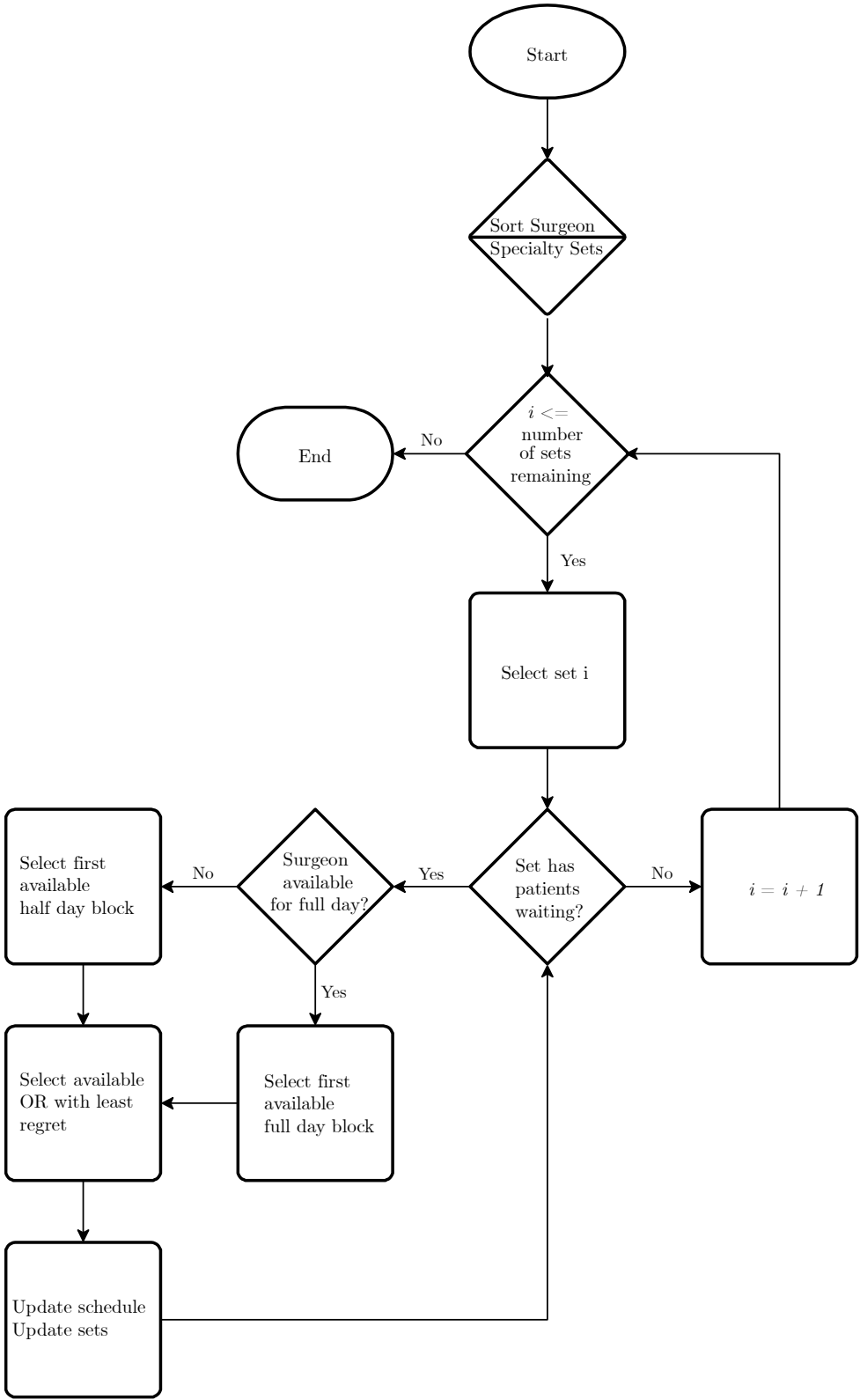

**Figure 2:** Constructive heuristic flow diagram.

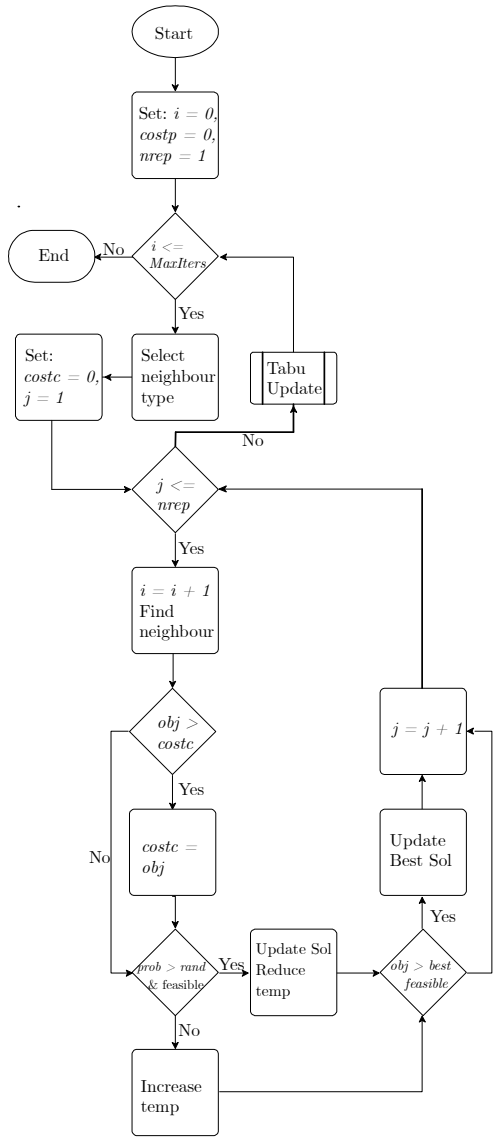 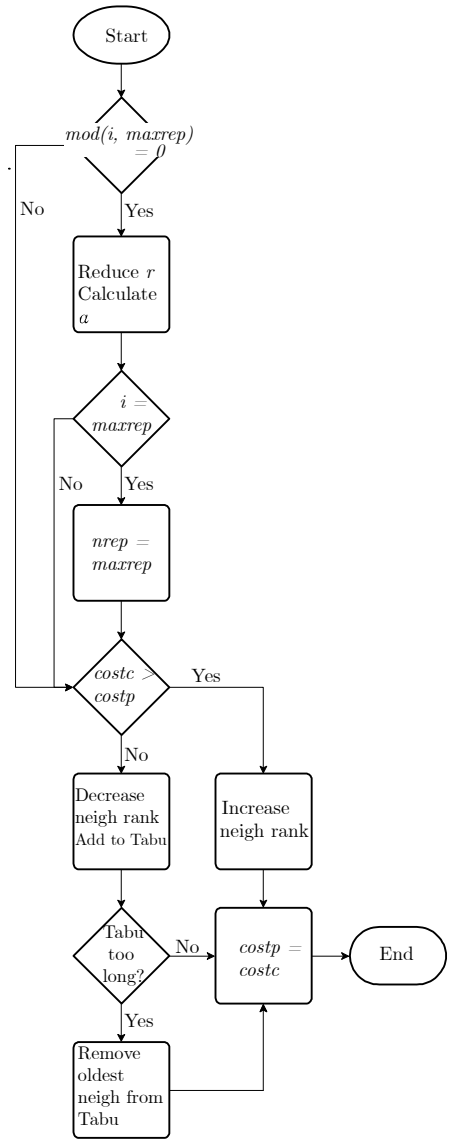

Figure 3: Hyper SA flow diagram.  Figure 4: Tabu Update flow diagram.

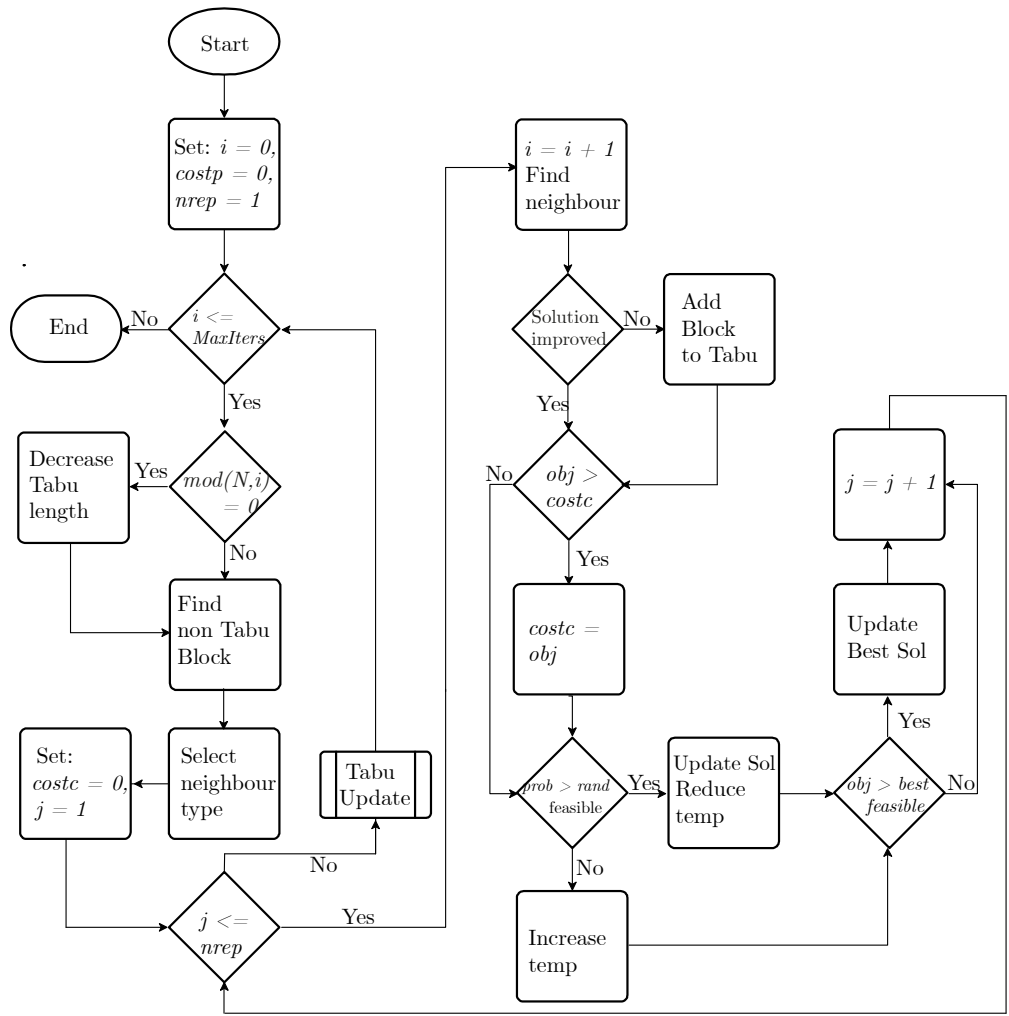

**Figure 5:** Hyper SA-TS Flow Diagram.